\documentclass[10pt]{amsart}

\usepackage[latin1]{inputenc}
\usepackage{amsmath}
\usepackage{amsfonts}
\usepackage{amssymb}
\usepackage[mathscr]{eucal}
\usepackage{diagrams}
\usepackage{xy,color,graphicx}
\usepackage{hyperref}
\xyoption{all}

\newcommand{\wis}[1]{{\text{ \usefont{OT1}{cmtt}{m}{n} #1}}}

\newcommand{\N}{\mathbb{N}}
\newcommand{\C}{\mathbb{C}}
\newcommand{\Z}{\mathbb{Z}}
\newcommand{\vtx}[1]{*+[o][F-]{\scriptscriptstyle #1}}
\newcommand{\Oscr}{\mathcal{O}}

\newtheorem{theorem}{Theorem}
\newtheorem{lemma}{Lemma}
\newtheorem{example}{Example}

\title{The singularities of noncommutative manifolds}
\author{Lieven Le Bruyn} 
\address{Department of Mathematics, University of Antwerp \\ 
 Middelheimlaan 1, B-2020 Antwerp (Belgium) \\ {\tt lieven.lebruyn@uantwerpen.be}}

\begin{document}
\sloppy

\maketitle

\begin{abstract}
We present a faster method to determine all singularities of quiver moduli spaces up to smooth equivalence. We show that every quiver controls a large family of noncommutative compact manifolds.
\end{abstract}

\section{The problem}

Let $Q$ be a quiver on a finite set of vertices $Q_v=\{ 1,2,\hdots,k \}$ having a finite set of arrows $Q_a$. The quiver may contain loops and oriented cycles. The structure of $Q$ is fully encoded in the Ringel bilinear form $\chi_Q$ on $\Z^{\oplus k}$ determined by
\[
\chi_Q(\epsilon_i,\epsilon_j) = \delta_{ij} - \# \{ a \in Q_a~|~\xymatrix{\circ_i \ar[r]^a & \circ_j} \} \]
The path algebra $\C Q$ has as $\C$-basis the set of all oriented paths in $Q$ including those of length zero which correspond to the vertices. Multiplication in $\C Q$ is induced by concatenation of paths. The path algebras $\C Q$ are very special cases of quasi-free, or formally smooth, algebras as in \cite{CuntzQuillen} or \cite{KontRos} and can be seen as corresponding to all possible noncommutative affine spaces among all noncommutative manifolds, see for example \cite{LBBook}. The study of the additive category of finite dimensional representations of $\C Q$ reduces to that of quiver-representations of $Q$. Each such representation has a dimension vector $\alpha=(a_1,\hdots,a_k) \in \N^{\oplus k}$ (giving the dimensions of the vertex-spaces) of total dimension $d(\alpha) = \sum_i a_i$. The set of all $Q$-representations of dimension vector $\alpha$ is an affine space $\wis{rep}_{\alpha}(Q)$ on which the group $GL(\alpha) = \prod_i GL_{a_i}$ acts via base-change in the vertex-spaces. It is well known that the corresponding GIT-quotient $\wis{rep}_{\alpha}(Q)/GL(\alpha)$ classifies isomorphism classes of semi-simple $Q$-representation of dimension vector $\alpha$ and that its coordinate ring, the ring of polynomial quiver invariants, is generated by traces along loops and oriented cycles in $Q$, see \cite{LBProcesi}.

Consider a stability structure $\theta=(\theta_1,\hdots,\theta_k) \in \Z^{\oplus k}$, then we call a $Q$-representation $V$ of dimension vector $\alpha$ a $\theta$-semistable representation if $\theta.\alpha = \sum_i \theta_i a_i = 0$ and if for every proper subrepresentation $W$ of $V$ we have $\theta.\beta_W \geq 0$ where $\beta_W$ is the dimension vector of $W$, and, if all $\theta.\beta_W > 0$ we say that $V$ is a $\theta$-stable representation. The corresponding moduli space $\wis{mod}_{\alpha}^{ss}(Q,\theta)$ of $\theta$-semistable representations of dimension $\alpha$ was introduced and studied in \cite{King}. Its points correspond to isomorphism classes of $\alpha$-dimensional representations of the form
\[
M = N_1^{\oplus e_1} \oplus \hdots \oplus N_u^{\oplus e_u} \]
such that all factors $N_i$ are $\theta$-stable of dimension vector $\beta_i$ and occur in $M$ with multiplicity $e_i \geq 0$, so that $\alpha = \sum_i e_i \beta_i$. It is well known  that most of these quiver moduli spaces $\wis{mod}_{\alpha}^{ss}(Q,\theta)$ are singular. In fact, in \cite{Domokos} it is shown that the only quivers $Q$ having all their quiver moduli spaces smooth are the Dynkin or extended Dynkin quivers.

Hence, we would like to determine the types of singularities that occur in $\wis{mod}_{\alpha}^{ss}(Q,\theta)$ up to smooth equivalence. Recall that points $x \in X$ and $y \in Y$ in two varieties $X$ and $Y$ are called smoothly equivalent, if there are natural numbers $k$ and $l$ and an isomorphism of complete local rings
\[
\widehat{\Oscr}_{X,x}[[x_1,\hdots,x_k]] \simeq \widehat{\Oscr}_{Y,y}[[y_1,\hdots,y_l]] \]
In principle, we can determine these finite number of types by combining the method of local quivers from \cite{Adri-LB} with Bocklandt's reduction steps \cite{Bocklandt} and \cite{SOS}.
In \cite{Adri-LB} the \'etale local structure of $\wis{mod}_{\alpha}^{ss}(Q,\theta)$ near $M$ was described as the quiver quotient variety of a local quiver setting $(Q_M,\alpha_M)$. Here, $Q_M$ is a quiver on $k$ vertices corresponding to the distinct stable factors of $M$. The number of directed arrows (or loops) from vertex $\circ_i$ to vertex $\circ_j$ is equal to $\delta_{ij}-\chi_Q(\beta_i,\beta_j)$. The main result of \cite{Adri-LB} asserts that there is an \'etale isomorphism between a neighborhood of $M$ in $\wis{mod}_{\alpha}^{ss}(Q,\theta)$ and a neighborhood of $\overline{0}$ in the quiver quotient variety $\wis{rep}_{\alpha}(Q_M)/GL(\alpha_M)$.
In \cite{SOS} Bocklandt's eduction steps from \cite{Bocklandt} were used to to classify such quiver quotient singularities up to smooth equivalence.

However, in all but the more trivial situations, this is a very time consuming method. Whereas Bocklandt's reduction steps are fairly efficient, the determination of all possible representation types of $M$ and the calculation of all local quiver settings is not. In this paper we introduce two concepts to speed up this process. First, we introduce an auxiliary quiver $Q_{\theta}$, depending only on the stability structure $\theta$ and not on the particular dimension vector $\alpha$, controlling all possible local quivers $(Q_M,\alpha_M)$. The quiver $Q_{\theta}$ will allow us to quickly determine the 'worst' singularity types in $\wis{mod}_{\alpha}^{ss}(Q,\theta)$. Next, we introduce a partial ordering on the possible types of quiver quotients, which can in any concrete situation be efficiently constructed inductively, and, which we use to characterize all other singularity types in $\wis{mod}_{\alpha}^{ss}(Q,\theta)$, starting from the 'worst' ones.

\section{The controlling quiver $Q_{\theta}$}

We fix a quiver $Q$ with vertices $Q_v= \{ 1,\hdots,k \}$ and fix a stability structure $\theta \in \N^{\oplus k}$. Let $\Sigma_{\theta}$ be the additive sub-monoid of $\N^{\oplus k}$ consisting of all dimension vectors $\alpha$ such that $\theta.\alpha = 0$. As the direct sum of two $\theta$-semistable representations is again $\theta$-semistable, we can consider in $\Sigma_{\theta}$ the additive sub-monoid $V_{\theta}$ consisting of those $\alpha \in \Sigma_{\theta}$ such that there exist $\theta$-semistable representations of $Q$ of dimension vector $\alpha$. In \cite{Schofield} an inductive procedure is given to determine $V_{\theta}$. Let $\{ \gamma_1,\hdots,\gamma_l \}$ be a minimal set of additive monoid generators of $V_{\theta}$. Such generating dimension vectors $\gamma$ have special properties:

\begin{lemma} With notations as above we have:
\begin{enumerate}
\item{$\gamma$ is in a minimal generator set of $V_{\theta}$ if and only if every $\theta$-semistable $Q$-representation is actually $\theta$-stable.}
\item{$\gamma$ is in a minimal generator set of $V_{\theta}$ if and only if the GIT-quotient
\[
\wis{rep}_{\alpha}^{ss}(Q) \rOnto \wis{mod}_{\alpha}^{ss}(Q,\theta) \]
is a principal $PGL_n$-fibration in the \'etale topology, where $\wis{rep}_{\alpha}^{ss}(Q)$ is the Zariski open subset of $\wis{rep}_{\alpha}(Q)$ consisting of $\theta$-semistable representations.}
\item{If $\gamma$ is in a minimal generator set of $V_{\theta}$ then the moduli space $\wis{mod}_{\alpha}^{ss}(Q,\theta)$ is smooth of dimension $1-\chi_Q(\gamma,\gamma)$.}
\end{enumerate}
\end{lemma}

\begin{proof}
Every $\theta$-semistable representation $M$ has a filtration by $\theta$-semistable sub-representations such that all filtration quotients $N_i = F{i+1}/F_i$ are $\theta$-stable. It follows that the $\theta$-semistable representation $N = \oplus_i N_i$ lies in the closure of the $GL(\alpha)$-orbit of $M$ and if the dimension vector of the $\theta$-stable factor $N_i$ is $\beta_i$ then clearly $\alpha = \sum_i \beta_i$. (1) follows from this, as does (2) using the fact that the stabilizer subgroup of a $\theta$-stable representation is $\C^*$. (3) follows from (2) and the fact that $\wis{rep}^{ss}_{\alpha}(Q)$ is a non-empty Zariski open subset of the affine space $\wis{rep}_{\alpha}(Q)$ is smooth.
\end{proof}

$Q_{\theta}$ will then be the quiver with vertices $\{ 1,\hdots,l \}$ where vertex $\circ_i$ corresponds to the generator $\gamma_i$. In $Q_{\theta}$ the number of directed arrows (or loops) from vertex $\circ_i$ to vertex $\circ_j$ will be equal to $\delta_{ij}-\chi_Q(\gamma_i,\gamma_j)$. A first use of $Q_{\theta}$ lies in the characterization of $\theta$-stable dimension vectors:

\begin{lemma} \label{lemma1} The following are equivalent
\begin{enumerate}
\item{There exists a $\theta$-stable representation of dimension vector $\alpha$}
\item{We can write $\alpha = \sum_{i=1}^l c_i \gamma_i$ where $\gamma=(c_1,\hdots,c_l)$ is the dimension vector of a simple representation of $Q_{\theta}$.}
\end{enumerate}
Moreover, in this case we have $\chi_Q(\alpha,\alpha) = \chi_{Q_{\theta}}(\gamma,\gamma)$.
\end{lemma}

\begin{proof} This is a direct consequence of \cite[Thm. 5.1]{Adri-LB} and the fact that $\{ \gamma_1,\hdots,\gamma_l \}$ generate $V_{\theta}$. The conclusion follows because the dimension of $\wis{mod}_{\alpha}^{ss}(Q,\theta)$ is in this case equal to $1-\chi_Q(\alpha,\alpha)$ and, on the other hand, the quiver setting $(Q_{\theta},\gamma)$ is the local quiver encoding the \'etale local structure of $\wis{mod}_{\alpha}^{ss}(Q,\theta)$ near a point
\[
N = N_1^{\oplus c_1} \oplus \hdots \oplus N_l^{\oplus l} \]
where $N_i$ is a $\theta$-stable representation of dimension vector $\gamma_i$. As $\gamma$ is a simple dimension vector for $Q_{\theta}$ the dimension of the quiver quotient variety $\wis{rep}_{\gamma}(Q_{\theta})/GL(\gamma)$ is equal to $1 - \chi_{Q_{\theta}}(\gamma,\gamma)$.
\end{proof}

The main purpose of the auxiliary quiver $Q_{\theta}$ is that it controls all local quiver settings $(Q_M,\alpha_M)$ describing the \'etale local structure of all moduli spaces $\wis{mod}_{\alpha}^{ss}(Q,\theta)$:

\begin{theorem} The quiver $Q_{\theta}$ contains enough information to construct the local quiver setting $(Q_M,\alpha_M)$ describing the \'etale local structure of the quiver moduli space $\wis{mod}_{\alpha}^{ss}(Q,\theta)$ near the point corresponding to
\[
M = S_1^{\oplus e_1} \oplus \hdots \oplus S_u^{\oplus e_u} \]
where the $S_i$ are non-isomorphic $\theta$-stable representations of dimension vector $\beta_i$. More precisely, if $N$ and $N'$ are $\theta$-stable representations of dimension vectors $\beta$ and $\beta'$ and if we can write $\beta=\sum_i c_i \gamma_i$ and $\beta' = \sum_i c'_i \gamma_i$ with $\gamma=(c_1,\hdots,c_l)$ and $\gamma'=(c'_1,\hdots,c'_l)$ simple dimension vectors of $Q_{\theta}$ then we have
\[
\chi_Q(\beta,\beta') = \chi_{Q_{\theta}}(\gamma,\gamma') \]
\end{theorem}

\begin{proof}
If $N \simeq N'$ (and hence $\beta=\beta'$ and $\gamma=\gamma'$) the claim follows from the previous lemma. So, assume that $N \not\simeq N'$, then the local quiver setting describing the \'etale local neighborhood of $\wis{mod}_{\beta+\beta'}^{ss}(Q,\theta)$ near $N \oplus N'$ is
\[
\xymatrix{\vtx{1} \ar@{=>}@(ul,dl)_{1-\chi_Q(\beta,\beta)} \ar@{=>}@/^2ex/[rr]^{-\chi_Q(\beta,\beta')} & & \vtx{1} \ar@{=>}@(ur,dr)^{1-\chi_Q(\beta',\beta')} \ar@{=>}@/^2ex/[ll]^{-\chi_Q(\beta',\beta)}} \]
On the other hand, $(Q_{\theta},\gamma+\gamma')$ is the local quiver describing the \'etale local structure of $\wis{mod}_{\beta+\beta'}^{ss}(Q,\theta)$ near $P \oplus Q$ where
\[
P = N_1^{\oplus c_1} \oplus \hdots \oplus N_l^{\oplus c_l} \quad \text{and} \quad Q = N_1^{\oplus c_1'} \oplus \hdots \oplus N_l^{\oplus c_l'} \]
By the previous lemma there are simple $Q_{\theta}$-representations $S$ and $T$ of dimension vector $\gamma$ and $\gamma'$ such that the semi-simple $Q_{\theta}$-representation $S \oplus T$ lies in a Zariski neighborhood of $\overline{0}$ in $\wis{rep}_{\gamma+\gamma'}(Q_{\theta})/GL(\gamma+\gamma')$ and hence corresponds to a point in $\wis{mod}_{\beta+\beta'}^{ss}(Q,\theta)$ corresponding to a representation $N_1 \oplus N_1'$ with $N_1$ and $N'_1$ both $\theta$-stable and of dimension vectors $\beta$ and $\beta'$. By the theory of local quivers of semi-simple quiver representations as in \cite{LBProcesi} the \'etale local structure of the quiver quotient variety near $S \oplus T$ is determined by the local quiver setting
\[
\xymatrix{\vtx{1} \ar@{=>}@(ul,dl)_{1-\chi_{Q_{\theta}}(\gamma,\gamma)} \ar@{=>}@/^2ex/[rr]^{-\chi_{Q_{\theta}}(\gamma,\gamma')} & & \vtx{1} \ar@{=>}@(ur,dr)^{1-\chi_{Q_{\theta}}(\gamma',\gamma')} \ar@{=>}@/^2ex/[ll]^{-\chi_{Q_{\theta}}(\gamma',\gamma)}} \]
As local quivers-settings only depend on the representation type, this quiver must be the same as the one of $N \oplus N'$ finishing the proof.
\end{proof}

From now on we will restrict attention to the study of moduli spaces $\wis{mod}_{\theta}^{ss}(Q,\theta)$ for dimension vectors $\alpha$ allowing $\theta$-stable representations. The general case reduces to this by the theory of general representations developed  in \cite{Schofield}.

\section{Bocklandt's reduction steps}

With $\wis{simps}$ we denote the set of all simple quiver settings, that is, all couples $(Q,\alpha)$ consisting of a quiver $Q$ and dimension vector $\alpha=(a_1,\hdots,a_k)$ with all $a_i \not= 0$ (that is, the support $supp(\alpha)$ of $\alpha$ contains all vertices of $Q$), which satisfy (see \cite{LBProcesi}):
\begin{itemize}
\item{$Q$ must be strongly connected, meaning that there exist directed paths connecting any two of its vertices, and,
\[
\chi_Q(\alpha,\epsilon_i) \leq 0 \quad \text{and} \quad \chi_Q(\epsilon_i,\alpha) \leq 0 \]
for all vertex dimension vectors $\epsilon_i$. That is, in every vertex $\circ_i$ the total number of incoming (and outgoing) dimensions is greater than or equal to the vertex-dimension.}
\item{If however $Q = \tilde{A}_k$ with cyclic orientation, then only $\alpha=(1,\hdots,1)$ is allowed.}
\end{itemize}
In verifying the numerical conditions it is practical to label each vertex with two numbers $\leq 0$,  giving the differences of the total incoming (resp outgoing) dimensions with the vertex-dimension. This allows us to spot quickly whether one of the reductions steps $(Q,\alpha) \rOnto (Q',\alpha')$, discovered by Raf Bocklandt in his characterization of smooth quiver quotient varieties, can be applied \cite{Bocklandt}:
\begin{itemize}
\item{If there is a vertex $\circ_i$ with vertex-dimension $a_i=1$ having loops, remove the loops to get $Q'$ and keep $\alpha'=\alpha$.}
\item{If one of the two numbers for $\circ_i$ is zero, remove the vertex $\circ_i$ and cable all arrows through, that is, any situation
\[
\xymatrix{\circ_k \ar@2[r]^a & \circ_i \ar@2[r]^b & \circ_l} \quad \text{becomes} \quad \xymatrix{\circ_k \ar@2[rr]^{a \times b} & & \circ_l} \]
to obtain $Q'$ and let $\alpha'$ be $\alpha$ with the $i$-th component removed.}
\item{If there is a vertex $\circ_i$ having a unique loop and such that one of the two numbers is $-1$, then we are in one of the following local situations in $\circ_i$
\[
\xymatrix{
& & \\
\vtx{1} \ar[r] & \vtx{k} \ar@(ul,ur) \ar@{.>}[ru] \ar@{.>}[rd] \ar@{.>}[r] & \\
& &} \qquad \xymatrix{
\ar@{.>}[rd] &  & \\
\ar@{.>}[r] & \vtx{k} \ar@(ul,ur)  \ar[r] & \vtx{1} \\
\ar@{.>}[ru] & & } \]
we replace the loop in $\circ_i$ by a bunch of $k$ arrows to or from $\xymatrix{\vtx{1}}$, that is, locally $Q'$ looks like
\[
\xymatrix{
& & \\
\vtx{1} \ar@{=>}[r]^k & \vtx{k}  \ar@{.>}[ru] \ar@{.>}[rd] \ar@{.>}[r] & \\
& &} \qquad \xymatrix{
\ar@{.>}[rd] &  & \\
\ar@{.>}[r] & \vtx{k}   \ar@{=>}[r]^k & \vtx{1} \\
\ar@{.>}[ru] & & }
\]
and keep $\alpha'=\alpha$.}
\end{itemize}
In every reduction step we either decrease the number of vertices or the number of loops. So, after a finite number of moves we arrive at a simple quiver setting $(Q^t,\alpha^t)$ which cannot be reduced further. As there is an element of choice in the reduction steps we can perform, there is a priori no reason that any two reduction procedures should result in the same final setting. Still, surprisingly, this is the case as was proved in \cite[\S 4]{SOS}. We will call this unique irreducible simple quiver setting $(Q^t,\alpha^t)$ the    type of $(Q,\alpha)$.

The upshot of this reduction process is the following result which follows from Bocklandt's work \cite{Bocklandt}:

\begin{theorem} \label{Zariski} For any $(Q,\alpha) \in \wis{simps}$ with (unique) type $(Q^t,\alpha^t)$ there is an isomorphism of varieties
\[
\wis{rep}_{\alpha}(Q)/GL(\alpha) \simeq \wis{rep}_{\alpha^t}(Q^t)/GL(\alpha^t) \times \C^d \]
where $d=\chi_{Q^t}(\alpha^t,\alpha^t)-\chi_Q(\alpha,\alpha)$. In particular, the corresponding quiver-quotient singularities are smoothly equivalent.
\end{theorem}

In \cite{Bocklandt} Bocklandt proves that the quiver quotient variety $\wis{rep}_{\alpha}(Q)/GL(\alpha)$ is smooth if and only if its type $(Q^t,\alpha^t)$ is either
\vskip 2mm
\[
\xymatrix{\vtx{1}}  \quad \text{or} \qquad  \xymatrix{\vtx{2} \ar@(dl,ul) \ar@(dr,ur)} \]
 
\vskip 5mm

\section{The partially ordered set of $\wis{types}$}

We will put a partial order on the set $\wis{types}$ of all types of simple quiver settings. Take an irreducible simple quiver setting $(Q,\alpha)$ and look for loops $I= \{ i \}$ or minimal oriented cycles $I = \{ i_1,\hdots, i_v \}$ in it. Consider the dimension vector $\beta_I = (\delta_{1I},\delta_{2I},\hdots,\delta_{kI})$ then there exists a simple $Q$-representation $S_I$ of dimension vector $\beta_I$. Now, consider the semi-simple $Q$-representation $M$ of dimension vector $\alpha = (a_1,\hdots,a_k)$
 \[
 M = S_I \oplus S_1^{\oplus a_1-\delta_{1I}} \oplus S_2^{\oplus a_2-\delta_{2I}} \oplus \hdots \oplus S_n^{\oplus a_n - \delta_{nI}} \]
 where $S_i$ is the $1$-dimensional simple vertex-representation in $\circ_i$. Let $(Q_M,\alpha_M)$ be the    local quiver setting associated to $M$ as defined in \cite{LBlocal}. In this case, $Q_M$ is a quiver on $k+1$ vertices $\{ 0,1,2,\hdots,k \}$, with $\circ_0$ corresponding to the simple component $S_I$, such that $Q_M | \{ 1,\hdots,k \} \simeq Q$ and the number of loops in $\circ_0$ is given by $1-\chi_Q(\beta_I,\beta_I)$ and the number of arrows from $\circ_0$ to $\circ_i$ (resp. from $\circ_i$ to $\circ_0$) is equal to $-\chi_Q(\beta_I,\epsilon_i)$ (resp. to $\chi_Q(\epsilon_i,\beta_I)$. The dimension vector $\alpha_M \in \N^{\oplus k+1}$ is determined by the multiplicities of the distinct simple factors in $M$, that is,
 \[
 \alpha_M = (1,a_1-\delta_{1I},\hdots,a_k-\delta_{kI}) \]
 Next, let $(Q'_M,\alpha'_M)$ be the (necessarily simple) quiver setting obtained by restricting $(Q_M,\alpha_M)$ to the support of $\alpha_M$. Finally, let $(Q_I,\alpha_I)$ be the type of $(Q'_M,\alpha'_M)$, then we say that $(Q_I,\alpha_I)$ is a direct successor of $(Q,\alpha)$ in $\wis{types}$ determined by the oriented cycle $I$ and we denote this by an arrow
\[
(Q,\alpha)  \rTo (Q_I,\alpha_I) \]
Composing such arrows then defines a partial order on $\wis{types}$.

\begin{example}
There is just one type of cycle (loop) $I = \{ 1 \}$ for the type 

\hskip 10mm 
\[
 \xymatrix{\vtx{2} \ar@(dl,ul) \ar@(dr,ur)} \]
 \hskip 5mm 
 
 \noindent 
 and the corresponding semi-simple representation $M$ is the direct sum $S_I \oplus S_1$ of two distinct simple $1$-dimensional representations, $S_I$ has one of the loops non-zero, $S_1$ not. The corresponding local quiver setting is then

\vskip 5mm
\[
\xymatrix{\vtx{1} \ar@(dl,d) \ar@(ul,u) \ar@/^2ex/[rr] & & \vtx{1} \ar@(u,ur) \ar@(d,dr) \ar@/^2ex/[ll]} \]

\vskip 10mm
\noindent
which has corresponding type $\xymatrix{\vtx{1}}$, that is, in $\wis{types}$ we have an arrow
\vskip 3mm
\[
\xymatrix{\vtx{2} \ar@(dl,ul) \ar@(dr,ur)} \hskip 10mm  \rTo \hskip 5mm  \xymatrix{\vtx{1}} \]

\vskip 3mm
\noindent
In fact, soon it will become apparent that $\xymatrix{\vtx{1}}$ is the unique minimal object in $\wis{types}$.
\end{example}

The upshot of this ordering is that it simplifies the singularity type as we move away from the worst singularity $\overline{0}$ in $\wis{rep}_{\alpha}(Q)/GL(\alpha)$ to singularities at points in the first deformed strata. 

Recall that points of the quiver quotient variety correspond to semi-simple representations of total dimension $\alpha$, that is, representations of the form
\[
N = T_1^{\oplus f_1} \oplus \hdots \oplus T_l^{\oplus f_l} \]
where all $T_i$ are simple $Q$-representations of dimension vector $\beta_i$ and such that $\sum_i f_i \beta_i = \alpha$. We then say that $N$ is of    representation type $\sigma(N) = (f_1,\beta_1;\hdots;f_l,\beta_l)$. The Luna stratification of $\wis{rep}_{\alpha}(Q)/GL(\alpha)$ consists of strata $\wis{strata}(\sigma)$, consisting of points of the same representation type $\sigma$, which are all locally closed subvarieties. In fact, one can show that $\wis{strata}(\sigma)$ is contained in the Zariski closure of $\wis{strata}(\sigma')$ if and only if the stabilizer subgroup $Stab(N_{\sigma'})$ is conjugated to a subgroup of $Stab(N_{\sigma})$in $GL(\alpha)$ for $N_{\sigma} \in \wis{strata}(\sigma)$ and $N_{\sigma'} \in \wis{strata}(\sigma')$.

The point $\overline{0}$ is contained in the most degenerate stratum of representation type $\sigma_0=(a_1,\epsilon_1;\hdots;a_n,\epsilon_n)$ if $\alpha = (a_1,\hdots,a_n)$. In the Hasse diagram of Luna strata, the strata of minimal dimension $\wis{strata}(\sigma)$ containing $\wis{strata}(\sigma_0)$ in their Zariski closure are exactly those with representation type of the form
\[
\sigma_I = (1,\beta_U;a_1-\delta_{1I},\epsilon_1;\hdots;a_n-\delta_{nI},\epsilon_n) \]
corresponding to a loop or a minimal proper oriented cycle $I$ in $Q$. The theory of local quivers, see for example \cite{LBlocal}, then asserts that the \'etale local structure of $\wis{rep}_{\alpha}(Q)/GL(\alpha)$ in a neighborhood of $M \in \wis{strata}(\sigma)$ is isomorphic to an \'etale local neighborhood of $\overline{0}$ in the quiver quotient variety $\wis{rep}_{\alpha_M}(Q_M)/GL(\alpha_M)$. 

Combining this with theorem~\ref{Zariski} we get the first assertion of the following result:

\begin{theorem} For $(Q,\alpha) \in \wis{types}$ we have:
\begin{enumerate} 
\item{If $(Q,\alpha) {\color{black} \rTo} (Q',\alpha')$ then any Zariski neighborhood  of $\overline{0}$ in the quotient variety $\wis{rep}_{\alpha}(Q)/GL(\alpha)$ contains points  smoothly equivalent with type $(Q',\alpha')$.}
\item{Every singularity of $\wis{rep}_{\alpha}(Q)/GL(\alpha)$ not of type $(Q,\alpha)$ is  smoothly equivalent to a singularity contained in $\wis{rep}_{\alpha'}(Q')/GL(\alpha')$ for some type $(Q',\alpha')$ such that $(Q,\alpha) \geq (Q',\alpha')$.}
\end{enumerate}
\end{theorem}

\begin{proof} As for the second assertion, recall that \'etale singularity types of quiver quotient varieties on it depend on their representation type and as $\overline{0}$ lies in the Zariski closure of any representation stratum, we have that any Zariski neighborhood of $\overline{0}$ contains points of all occurring \'etale singularity types. Now, take such a singularity type $\tau$ with associated representation type $\sigma$ and consider a representation type $\sigma'$ having a stratum of minimal dimension such that
\[
\wis{strata}(\sigma_0) \subsetneq \wis{strata}(\sigma') \subset \overline{\wis{strata}(\sigma)} \]
then $\sigma' = \sigma_I$ for some loop or minimal oriented cycle $I$ in $Q$. Any point in $\wis{strata}(\sigma)$ is of type $\tau$ and as $\wis{strata}(\sigma_I)$ is contained in its Zariski closure, the theory of local quivers entails that any Zariski neighborhood of $\overline{0}$ in $\wis{rep}_{\alpha_M}(Q_M)/GL(\alpha_M)$, for $M$ of type $\sigma_I$, contains points \'etale of type $\tau$. As  $\wis{rep}_{\alpha_M}(Q_M)/GL(\alpha_M)$ and $\wis{rep}_{\alpha'}(Q')/GL(\alpha')$ are smoothly isomorphic there are singularities in $\wis{rep}_{\alpha'}(Q')$ \'etale smoothly equivalent to $\tau$, finishing the proof.
\end{proof}
This then gives an  algorithm to determine all singularity types of quiver quotient varieties up to  smooth equivalence.

\begin{theorem} Let $(Q,\alpha) \in \wis{simps}$ and apply reduction steps to determine is type $(Q^t,\alpha^t) \in \wis{types}$. Then, the singularity types of points in the quiver quotient variety $\wis{rep}_{\alpha}(Q)/GL(\alpha)$ are, up to smooth equivalence, exactly those $(Q',\alpha') \in \wis{types}$ such that $(Q^t,\alpha^t) {\color{black} \geq} (Q',\alpha')$.
\end{theorem}

Note that different types may still be smoothly equivalent. For example, in \cite{SOS} we showed that types $5_{3a}$ and type $5_{4c}$ have isomorphic rings of invariants. Further, the locus of all points in $\wis{rep}_{\alpha}^{ss}(Q)/GL(\alpha)$ consisting of points of a specific type may consist of several strata $\wis{strata}(\sigma)$, even of varying dimensions.

\section{Hitchhiker's guide to $\wis{types}$}

In principle one can build a map of the partial ordered set $\wis{types}$, inductively by dimension of the quotient variety, and by number of the vertices in the quiver. If the dimension is $D$ and the number of vertices is $n$ we will enumerate all possible types $(Q,\alpha)$ as $D_{na},D_{nb},\hdots$. The quiver $Q$ is determined by the integral matrix $M_Q$  describing its Ringel bilinear form $\chi_Q$, and the condition that $(Q,\alpha) \in \wis{types}$ of dimension $D$ can then be expressed as a system of equations and inequalities involving the entries of $M_Q$ and $\alpha$, starting with
\[
1-\chi_Q(\alpha,\alpha)=D \qquad \chi_Q(\alpha,\epsilon_i) \leq 0 \qquad \chi_Q(\epsilon_i,\alpha) \leq 0 \]
followed by relations expressing that none of the Bocklandt's reduction steps are possible for $(Q,\alpha)$. This then produces a list of all types of dimension $D$ and we have to relate them to the already constructed poset of types of dimension $< D$.

This involves computing local quiver settings for representation types $\sigma_I$ corresponding to loops or minimal oriented cycles $I$ in the quiver $Q$. As there will be at least one loop in this local quiver in the vertex $\circ_0$, having vertex-dimension $1$, we see that by going to its associated type we drop he dimension of the quotient-variety by at least one. That is, we will only have to draw red arrows connecting the new types of dimension $D$ to types already constructed before.

We can also describe easily, for all dimensions $D$, the types having only an arrow to the unique minimal element $\xymatrix{\vtx{1}}$. However, we will not draw this arrow so that these types become sinks in the map. This happens when the corresponding quotient variety is an isolated singularity and those quiver-quotients have been classified in \cite{BLS} to be of the form
\[
\xymatrix@=.4cm{
\vtx{1} \ar@{.>}@/_6ex/[ddd]& \vtx{1} \ar@{=>}[l]^{k_4} & \\
& & \vtx{1} \ar@{=>}[lu]^{k_3} \\
& & \vtx{1} \ar@{=>}[u]^{k_2} \\
\vtx{1} \ar@{=>}[r]^{k_l} & \vtx{1} \ar@{=>}[ru]^{k_1} &} \]
where $Q$ has $l$ vertices and all $k_i \geq 2$. The resulting dimension is then $D=\sum_i k_i + l -1$. In \cite{SOS} all types of dimension $D \leq 6$ have been classified. The first dimension allowing a quiver-quotient singularity is $D=3$ and there is just one such type, corresponding to the conifold singularity $3_c$
\vskip 2mm
\[
\xymatrix{\vtx{1} \ar@{=>}@/^2ex/[rr] & & \vtx{1} \ar@{=>}@/^2ex/[ll]} \]
\vskip 3mm
\noindent
In dimension $D=4$ there are exactly three types
\vskip 3mm
\[ 4_2~:~ \xymatrix{ \vtx{1} \ar@{=>}@/^2ex/[rr] & & \vtx{1} \ar@3@/^2ex/[ll]}  \quad
4_{3a}~:~\xymatrix{\vtx{1}\ar@/^/[rr]\ar@/^/[rd]&&\vtx{1}\ar@/^/[ll]\ar@/^/[ld]\\
&\vtx{1}\ar@/^/[ru]\ar@/^/[lu]&} \quad 4_{3b}~:~\xymatrix{\vtx{1}\ar@2@/^/[rr]&&\vtx{1}\ar@2@/^/[ld]\\
&\vtx{1}\ar@2@/^/[lu]&} \]
\vskip 2mm

\begin{figure} \label{map}
\begin{center}
\[
\includegraphics[width=\textwidth]{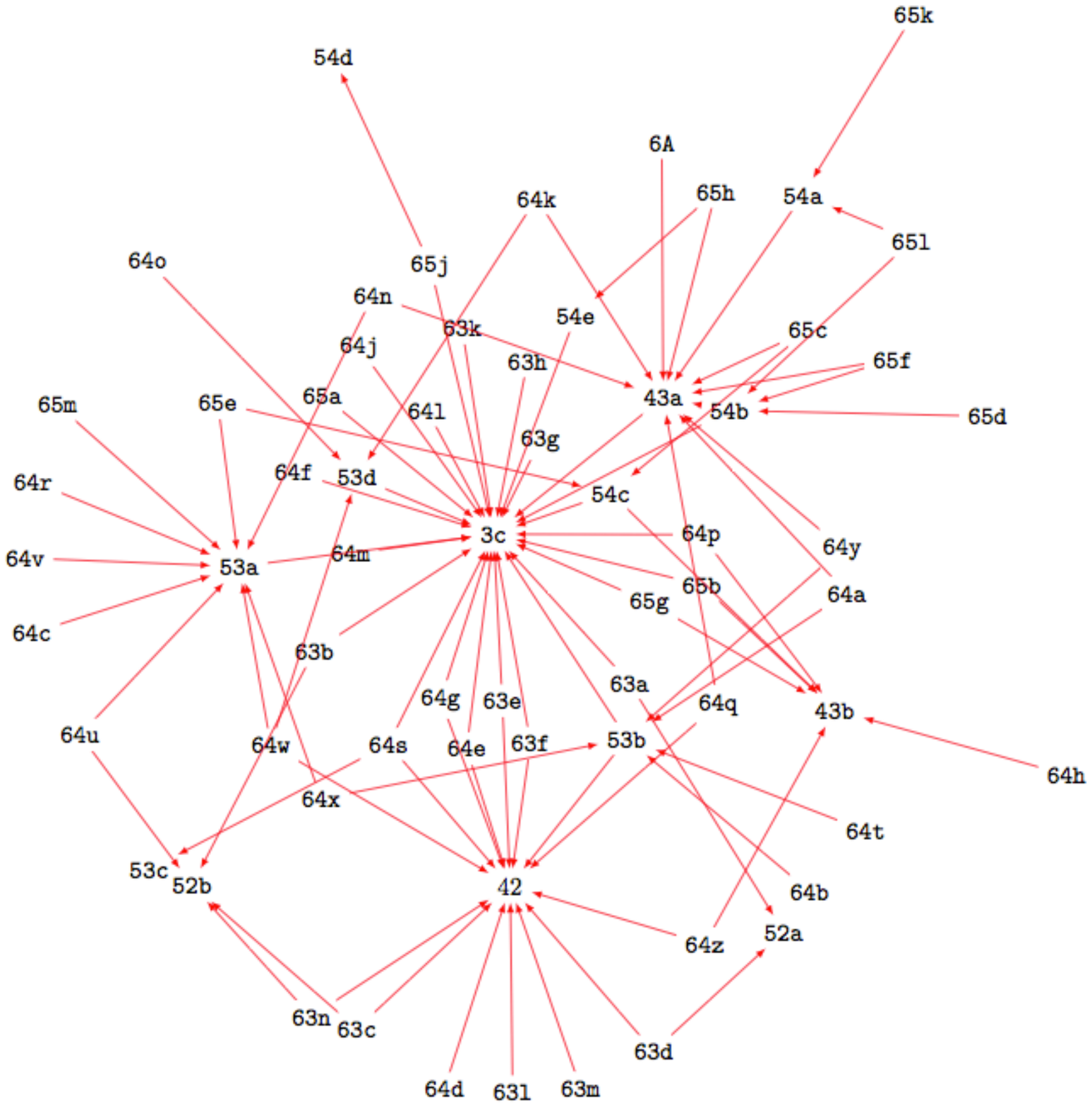}
\]
\end{center}
\caption{Hitchhiker's guide to  $\wis{types}$}
\end{figure}

\noindent
Dimension $D=5$ adds 11 types and in dimension $D=6$ we get an additional 54 types, see \cite{SOS} for all details. In \cite{SOS} these types where then classified up to isomorphism by a method called 'fingerprinting' singularities of which our partial order on types is a scaled down version. Using the enumeration of types as given in \cite{SOS} we can then draw the poset $\wis{types}$ up to dimension $D=6$ as in Figure~1. The isolated quiver-singularities in dimension $6$ are not included.

\section{Singularities of quiver moduli spaces}

We are now in a position to describe more efficiently the types of singularities that occur in $\wis{mod}_{\alpha}^{ss}(Q,\theta)$ up to smooth equivalence:

\begin{itemize}
\item{For fixed dimension vector $\alpha$ there is a limited number of possible simple dimension vectors $\{ \gamma_1,\hdots,\gamma_v \}$ of the controlling quiver $Q_{\theta}$ as in lemma~\ref{lemma1}.}
\item{For $1 \leq i \leq v$ apply Bocklandt's reduction steps to obtain the type $(Q_{\theta}^t,\gamma_i^t)$ of $(Q_{\theta} | supp(\gamma_i),\gamma_i)$.}
\item{Up to smooth equivalence, the types of singularities that occur in $\wis{mod}_{\alpha}^{ss}(Q,\theta)$ are precisely those $(Q',\alpha')$ such that
\[
(Q_{\theta}^t,\gamma_i^t) \geq (Q',\alpha') \]
for at least one $1 \leq i \leq v$.}
\end{itemize}

In particular, this allows us to characterize for a given quiver $Q$ and stability structure $\theta$ all the moduli spaces $\wis{mod}_{\alpha}^{ss}(Q,\theta)$ which are smooth, as those for which all types $(Q^t_{\theta},\gamma_i^t)$ are either 
\[ \xymatrix{\vtx{1}}~\qquad \text{or} \hskip 10mm \qquad~\xymatrix{\vtx{2} \ar@(dl,ul) \ar@(dr,ur)} \]
We will illustrate the above by some examples. First, we will determine the controlling quiver $Q_{\theta}$ relevant in the study of representations of the modular group. Then we will illustrate how one can inductively extend on Figure~1, and, finally we will give a short proof of the classification of all smooth quiver moduli spaces in this case.

\begin{example} \label{quiver}
Consider the quiver $Q$ below and stability structure $\theta=(-1,-1;1,1,1)$ that appears naturally in the study of character varieties of the modular group $\Gamma = PSL_2(\Z)$, see \cite{LBbraid}. 
\[
\xymatrix@=.25cm{
& & & & \vtx{x} \\
\vtx{a} \ar[rrrru] \ar[rrrrd] \ar[rrrrddd] & & & & \\
& & & & \vtx{y} \\
\vtx{b} \ar[rrrru] \ar[rrrruuu] \ar[rrrrd] & & & & \\
& & & &  \vtx{z} }
\]
then simple $\Gamma$-representations of dimension $n=a+b=x+y+z$ determine $\theta$-stable representations of dimension vector $\alpha=(a,b;x,y,z)$, which then must satisfy $min(a,b) \geq max(x,y,z)$. In this case, the monoid $V_{\theta}$ is generated by the six dimension vectors 
\[
\begin{cases} \gamma_1=(1,0;1,0,0) \quad \gamma_2=(0,1;0,1,0) \quad \gamma_3=(1,0;0,0,1) \\
\gamma_4=(0,1;1,0,0) \quad \gamma_5=(1,0;0,1,0) \quad \gamma_6=(0,1;0,0,1)
\end{cases} \]
which correspond to the six one-dimensional representations of $\Gamma = C_2 \ast C_3$. In this case the quiver $Q_{\theta}$ is 
\[
\xymatrix@=.5cm{
& \vtx{a_1} \ar@/^/[ld] \ar@/^/[rd] & \\
\vtx{a_6} \ar@/^/[ru]  \ar@/^/[d] & & \vtx{a_2} \ar@/^/[lu] \ar@/^/[d] \\
\vtx{a_5} \ar@/^/[u]  \ar@/^/[rd] & & \vtx{a_3} \ar@/^/[u] \ar@/^/[ld] \\
& \vtx{a_4} \ar@/^/[lu] \ar@/^/[ru]  & }
\]
and a dimension vector $\alpha_{\theta}=(a_1,\hdots,a_6)$ corresponds to $\alpha$ if and only if
\[
\begin{cases}
a=a_1+a_3+a_5 \quad b=a_2+a_4+a_6 \\
x=a_1+a_4 \quad y= a_2+a_5 \quad z=a_3+a_6
\end{cases}
\]
\end{example}

\begin{example} \label{many}
We will determine the singularity types occurring in the (unique) singular moduli space of smallest possible dimension, which is $7$, corresponding to dimension vector $\alpha=(3,3;2,2,2)$. Up to hexagonal symmetry there are two corresponding simple dimension vectors $\alpha_{\theta}^{(1)},\alpha_{\theta}^{(2)}$
\[
\xymatrix@=.5cm{
& \vtx{1} \ar@/^/[ld] \ar@/^/[rd] & \\
\vtx{1} \ar@/^/[ru]  \ar@/^/[d] & & \vtx{1} \ar@/^/[lu] \ar@/^/[d] \\
\vtx{1} \ar@/^/[u]  \ar@/^/[rd] & & \vtx{1} \ar@/^/[u] \ar@/^/[ld] \\
& \vtx{1} \ar@/^/[lu] \ar@/^/[ru]  & } \qquad \qquad \xymatrix@=.5cm{
& \vtx{2} \ar@/^/[ld] \ar@/^/[rd] & \\
\vtx{1} \ar@/^/[ru]  \ar@/^/[d] & & \vtx{2} \ar@/^/[lu] \ar@/^/[d] \\
\vtx{0} \ar@/^/[u]  \ar@/^/[rd] & & \vtx{1} \ar@/^/[u] \ar@/^/[ld] \\
& \vtx{0} \ar@/^/[lu] \ar@/^/[ru]  & }
\]
Neither of these quiver settings can be further reduced so they determine two new types, let us call them resp. $7_{6a}$ and $7_{4a}$ (recall that the first subindex gives the number of vertices of the quiver). Next, we have to connect them to the map of Figure~1.

We have to determine the local quivers associated to minimal oriented cycles in these two quivers. For type $7_{6a}$ there is up to symmetry just one such cycle, namely between two consecutive vertices. For type $7_{4a}$ we have up to symmetry two possible cycles, either containing an ending vertex or between the two middle vertices. 

The corresponding local quiver for type $7_{6a}$ and the second possibility for type $7_{4a}$ are both of the form
\vskip 3mm
\[
\xymatrix@=.4cm{& \vtx{1} \ar@(ul,ur) \ar@/^1ex/[ld] \ar@/^1ex/[rd] & \\
\vtx{1} \ar@/^1ex/[ru] \ar@/^1ex/[d] & & \vtx{1} \ar@/^1ex/[lu] \ar@/^1ex/[d] \\
\vtx{1} \ar@/^1ex/[u] \ar@/^1ex/[rr] & & \vtx{1} \ar@/^1ex/[ll] \ar@/^1ex/[u]} \]
and after reducing the loop we obtain type $6_{5k}$ of the classification from \cite{SOS}. The first possibility for type $7_{4a}$ gives as local quiver
\vskip 3mm
\[
\xymatrix{& \vtx{1} \ar@(ul,ur) \ar@/^1ex/[d] & \\
\vtx{1} \ar@/^1ex/[r] & \vtx{2} \ar@/^1ex/[l] \ar@/^1ex/[r] \ar@/^1ex/[u] & \vtx{1} \ar@/^1ex/[l]} \]
and, after reducing the loop, we obtain type $6_A$ from \cite{SOS}. 

Both type $6_{5k}$ and $6_A$ are already on our map (upper right hand), so we obtain that up to smooth equivalence the singularities of $\wis{mod}_{\alpha}^{ss}(Q,\theta)$ are exactly of the types
\[
7_{6a},~7_{4a},~6_{5k},~6_A,~5_{4a},~4_{3a},~3_{con} \]
with degeneration diagram
\[
\xymatrix{7_{6a} {\color{red} \ar[r]} & 6_{5k} {\color{red} \ar[r]} & 5_{4a} {\color{red} \ar[rd]} & & \\
7_{4a} {\color{red} \ar[r] \ar[ru]} & 6_A {\color{red} \ar[rr]} & & 4_{3a} {\color{red} \ar[r]} & 3_{con}}
\]
This calculation illustrates the iterative process involved to extend on Figure~1 in a specific problem.
\end{example}

\begin{example}
Next, let us see how to apply the foregoing in order to classify the smooth quiver moduli spaces. Before, we have seen that for $\alpha=(3,3;2,2,2)$ in the moduli space there should be singularities smoothly equivalent to the conifold singularity. Indeed, if we take representations of the form $M = N \oplus N'$ where $N$ (resp. $N'$) is $\theta$-stable of dimension vector $(2,1;1,1,1)$ (resp. $(1,2;1,1,1)$) then the corresponding local quiver is
\vskip 4mm
\[
\xymatrix{\vtx{1} \ar@(u,ul) \ar@(d,dl) \ar@2@/^2ex/[rr] & & \vtx{1} \ar@(u,ur) \ar@(d,dr) \ar@2@/^2ex/[ll]}
\]
\vskip 7mm
\noindent
Now, if $\beta$ is a strictly larger dimension-vector (meaning that all its vertex dimensions are greater or equal than those of  $\alpha$) we can write
\[
\beta = \alpha + \sum_{i=1}^6 n_i \gamma_i \]
and calculating the local quiver of this representation type shows that we can never reduce to a smooth setting. Hence, all moduli spaces for dimension vectors larger than $\alpha$ will be singular. The remaining dimension vectors are either Dynkin or extended Dynkin (and hence have smooth moduli spaces) or of the form
\[
(b+1,b;b,b,1), \quad (b,b;b,b-1,1) \quad \text{or} \quad (4,2;2,2,2) \]
In the first case, there are (up to symmetry) two possible families of simple dimension vectors for $Q_{\theta}$, namely
\[
\xymatrix@=.5cm{
& \vtx{x} \ar@/^/[ld] \ar@/^/[rd] & \\
\vtx{0} \ar@/^/[ru]  \ar@/^/[d] & & \vtx{x} \ar@/^/[lu] \ar@/^/[d] \\
\vtx{b-x} \ar@/^/[u]  \ar@/^/[rd] & & \vtx{1} \ar@/^/[u] \ar@/^/[ld] \\
& \vtx{b-x} \ar@/^/[lu] \ar@/^/[ru]  & } \qquad \qquad \xymatrix@=.5cm{
& \vtx{b-x} \ar@/^/[ld] \ar@/^/[rd] & \\
\vtx{1} \ar@/^/[ru]  \ar@/^/[d] & & \vtx{b-1-x} \ar@/^/[lu] \ar@/^/[d] \\
\vtx{1+x} \ar@/^/[u]  \ar@/^/[rd] & & \vtx{0} \ar@/^/[u] \ar@/^/[ld] \\
& \vtx{x} \ar@/^/[lu] \ar@/^/[ru]  & }
\]
The first one we can reduce to
\vskip 2mm
\[
\xymatrix{\vtx{x} \ar@(ul,dl) \ar@/^2ex/[r] & \vtx{1} \ar@/^2ex/[l] \ar@/^2ex/[r] & \vtx{b-x} \ar@(ur,dr) \ar@/^2ex/[l]} \]
\vskip 2mm
\noindent
and subsequently to
\vskip 2mm
\[
\xymatrix{\vtx{x} \ar@/^2ex/[r] & \vtx{1} \ar@{=>}@/^2ex/[l]^x \ar@{=>}@/^2ex/[r]^{b-x} & \vtx{b-x} \ar@/^2ex/[l]} \]
\vskip 2mm
\noindent
and finally to
\vskip 2mm
\[
\xymatrix{\vtx{1} \ar@(ul,dl)_x \ar@(ur,dr)^{b-x} } \]
\vskip 3mm
\noindent
which is of type $\xymatrix{\vtx{1}}$ whence smooth. The second one first reduces to
\vskip 5mm
\[
\xymatrix{\vtx{x} \ar@(ul,dl) \ar@/^2ex/[r] & \vtx{1} \ar@(ul,ur) \ar@(dl,dr) \ar@/^2ex/[l] \ar@/^2ex/[r] & \vtx{b-1-x} \ar@(ur,dr) \ar@/^2ex/[l]} \]
\vskip 5mm
\noindent
which after reducing the two loops in the middle vertex is of the same type as the first reduction of the first case, so is again smooth. For the dimension vectors $(b,b;b,b-1,1)$ the argument is similar. As for the special dimension vector $\alpha=(4,2;2,2,2)$ here the corresponding quiver-setting $(Q_{\theta},\alpha_{\theta})$ is
\vskip 3mm
\[
\xymatrix{\vtx{1} \ar@/^2ex/[r] & \vtx{1} \ar@/^2ex/[l] \ar@/^2ex/[r] & \vtx{2} \ar@/^2ex/[l] \ar@/^2ex/[r] & \vtx{1} \ar@/^2ex/[l] \ar@/^2ex/[r] & \vtx{1} \ar@/^2ex/[l]} \]
\vskip 3mm
\noindent
which is easily seen to reduce to $\xymatrix{\vtx{1}}$, whence is smooth. 
\end{example}

Concluding, we have the following characterization of all smooth moduli spaces for the above quiver $Q$ and stability structure $\theta=(-1,-1;1,1,1)$, see also \cite{ABV}:

\begin{theorem} $\wis{mod}_{\alpha}^{ss}(Q,\theta)$ is smooth unless all vertex-dimensions of $\beta$ are greater or equal than those of $\beta=(3,3;2,2,2)$.
\end{theorem}

Observe from \cite{LBbraidreversion} that these are exactly the components on which transposition induces the identity.

\section{Compact noncommutative manifolds}

If the quiver $Q$ has no oriented cycles all moduli spaces $\wis{mod}_{\theta}^{ss}(Q,\theta)$ are projective varieties and in \cite{LBNCM} it is argued that one can view the family of projective varieties
\[
(~\bigsqcup_{d(\alpha)=n} \wis{mod}_{\alpha}^{ss}(Q,\theta)~)_n \]
as a noncommutative compact manifold. That is, the additive category $\wis{rep}^{ss}(Q,\theta)$ of all finite dimensional $\theta$-stable representations of $Q$ can be covered by representation categories $\wis{rep}(A)$ consisting of all finite dimensional representations of formally smooth algebras $A$. We have seen that the singularities of the noncommutative compact manifold, as well as all local quiver settings describing its \'etale local quiver, are controlled by the quiver $Q_{\theta}$. Conversely, we have

\begin{theorem} For every quiver $Q^{\dagger}$ there exist noncommutative compact manifolds of the form
\[
(~\bigsqcup_{d(\alpha)=n} \wis{mod}_{\alpha}^{ss}(Q,\theta)~)_n \]
with $Q$ having no oriented cycles such that all local quiver settings are controlled by $Q^{\dagger}$.
\end{theorem}

\begin{proof}
We start with $Q_0=Q^{\dagger}$ and stability structure $\theta_0=(0,\hdots,0)$. We use the trick of iterating the procedure of doubling vertices, see \cite[\S 2]{Domokos}, to remove all loops and oriented cycles in and to modify the stability structure accordingly. That is, if after $k$ steps we have arrived at a situation $(Q_k,\theta_k)$ such that all local quiver settings for moduli spaces of $\theta_k$-semistable representations are controlled by $Q^{\dagger}$ and if we still have a vertex $\circ_i$ in $Q_k$ having loops or oriented cycles passing through it, and if the $i$-th $\theta_k$-component is $t_i$, then we modify the situation by splitting the vertex in two vertices $\circ_{i_-}$ and $\circ_{i_+}$ and adjusting loops and arrows starting or ending in $\circ_i$ as indicated below
\[
\xymatrix{ \ar@{.>}[rd] & &  \\
 \ar@{.>}[r]& \circ_i \ar@{.>}[ru] \ar@{.>}[r] \ar@{.>}[rd] \ar@{.>}@(ul,ur) &  \\
 \ar@{.>}[ru] & & } \qquad~\qquad
 \xymatrix{ \ar@{.>}[rd] & & &  \\
 \ar@{.>}[r]& \circ_{i_+} & \circ_{i_-} \ar@{.>}[ru] \ar@{.>}[r] \ar@{.>}[rd] \ar[l] \ar@{.>}@/_2ex/[l] &  \\
 \ar@{.>}[ru] & & & }
\]
to get a new quiver $Q_{k+1}$ and new stability structure $\theta_{k+1}$ which coincides with $\theta_k$ in all non-modified vertices and is equal to $-n$ in $\circ_{i_-}$ and equal to $t_i+n$ in $\circ_{i_+}$, where $n$ is chosen large enough to ensure that all local quiver settings for moduli spaces of $\theta_{k+1}$-semistable representations of $Q_{k+1}$ are controlled by $Q^{\dagger}$. Note that if we have a dimension vector $\alpha_k$ allowing $\theta_k$-semistable representation, then the dimension vector $\alpha_{k+1}$, which coincides with $\alpha_k$ in all non-modified vertices and with the $\alpha_k$-component in $\circ_i$ in the new vertices $\circ_{i_-}$ and $\circ_{i_+}$, will allow $\theta_{k+1}$-semistable representations. By \cite[Thm. 2.2]{Domokos} this is always possible. In fact, we can even choose $n$ such that the local quiver settings for $\wis{mod}_{\alpha_k}^{ss}(Q_k,\theta_k)$ are exactly  local quiver settings  for $\wis{mod}_{\alpha_{k+1}}^{ss}(Q_{k+1},\theta_{k+1})$. After a finite number of steps we obtain a quiver $Q$ having no loops nor oriented cycles and a stability structure $\theta$ such that the corresponding noncommutative compact manifold is controlled by $Q^{\dagger}$.
\end{proof}

As is clear from the foregoing proof, the same quiver $Q^{\dagger}$ can control a large array of compact noncommutative manifolds, which can be quite different.

\begin{example}
Take as quiver $Q^{\dagger} = Q_{\theta}$ of example~\ref{quiver} controlling the noncommutative compactification of the modular group. This quiver also controls the noncommutative compact manifold defined by the quiver $Q'$ below
\[
\xymatrix{& & \vtx{} \ar[ld] \ar[rd] & & \\
\vtx{} \ar[r] \ar[rru] \ar[d] & \vtx{} & & \vtx{} & \vtx{} \ar[llu] \ar[d] \ar[l] \\
\vtx{} \ar[ru] \ar[rr] & & \vtx{} & & \vtx{} \ar[lu] \ar[ll] \\
& & \vtx{} \ar[llu] \ar[rru] \ar[u] & &} 
\]
and stability structure (with cyclic ordering op vertices of $Q^{\dagger}$ and split vertices as consecutive entries)
\[
\theta' = (0;-p,p;0;-q,q;0;-r,r) \]
where $p,q$ and $r$ are sufficiently large primes. To a simple dimension vector $\alpha_{\dagger}=(a_1,a_2,a_3,a_4,a_5,a_6)$ of $Q^{\dagger}$ there is a unique dimension vector 
\[
\alpha' = (a_1;a_2,a_2;a_3;a_4,a_4;a_5;a_6,a_6) \]
of $Q'$ allowing $\theta'$-stable representations. By \cite[Thm. 2.2]{Domokos} the local quiver settings for the moduli space $\wis{mod}_{\alpha'}^{ss}(Q',\theta')$ are exactly the same as those of the quiver quotient-variety $\wis{rep}_{\alpha_{\dagger}}(Q^{\dagger})/GL(\alpha_{\dagger})$.

On the other hand, for the quiver $Q$ of example~\ref{quiver} and stability structure $\theta=(-1,-1;1,1,1)$ the local quivers for the moduli space $\wis{mod}_{\alpha}^{ss}(Q,\theta)$ for a dimension vector $\alpha$ allowing $\theta$-stable representations are, in general, determined by those of several simple dimension vectors $\alpha_{\dagger}$ of $Q^{\dagger}$ as example~\ref{many} illustrates.
\end{example}

\end{document}